\documentclass[12pt,twosided,reqno]{amsart}
\usepackage{amsmath}
\usepackage{amsthm}
\usepackage{amsfonts}
\usepackage{amssymb}

\theoremstyle{theorem}

\theoremstyle{definition}

\begin{document}

\vskip .5cm
\title{Area of  triangles associated with a curve}
 \vskip0.5cm
\thanks{
    \newline\indent { This was   supported by Basic Science Research Program through
    the National Research Foundation of Korea (NRF) funded by the Ministry of Education, Science and Technology (2010-0022926).}
}

\vskip 0.5cm

\maketitle

\vskip 0.5cm
\centerline{\scshape
Dong-Soo Kim and Kyu-Chul Shim} \vskip .2in

\begin{abstract}
It is well known that the area $U$ of the triangle
formed by three tangents to a parabola $X$ is half of the area $T$
 of the
triangle formed by joining their points of contact.
In this article, we study some properties
of $U$ and $T$ for strictly convex plane curves. As a result, we
establish a characterization for parabolas.

\end{abstract}

\vskip 1cm

\date{}
\maketitle

\section{Introduction}

 \vskip 0.50cm
Let $X=X(s)$ be a  unit speed smooth  curve in the  plane
 ${\mathbb R}^{2}$ with nonvanishing curvature,
and let $A =
X(s)$,$A_i=X(s+h_i), i = 1, 2,$ be three distinct neighboring points on $X$.
Denote by $\ell$, $\ell_1$, $\ell_2$ the tangent
lines passing through the points $A,A_1,A_2$ and by $B,B_1,B_2$
the intersection points $\ell_1\cap \ell_2$, $\ell\cap \ell_1$,
$\ell\cap \ell_2$, respectively.
It is well known that the area $U(s,h_1,h_2)=|\bigtriangleup BB_1B_2|$ of the triangle
formed by three tangents to a parabola is half of the area
$T(s,h_1,h_2)=|\bigtriangleup AA_1A_2|$ of the
triangle formed by joining their points of contact ([1]).

 \vskip 0.3cm
 The present article studies whether this property exhaustively characterizes parabolas.

Usually, a regular plane curve $X:I \rightarrow {\mathbb R}^{2}$  defined on an open interval
 is called {\it convex} if, for all $t\in I$, the trace $X(I)$ lies
 entirely on one side of the closed half-plane determined by the tangent line at $X(t)$ ([2]).

Hereafter,  we will  say  that a simple convex curve $X$ in the  plane
 ${\mathbb R}^{2}$ is  {\it strictly convex} if the curve   is smooth (that is, of class $C^{(3)}$) and is of positive  curvature $\kappa$
 with respect to the unit normal $N$ pointing to the convex side.
 Hence, in this case we have $\kappa(s)=\left< X''(s), N(X(s))\right> >0$, where $X(s)$ is an arclength parametrization of $X$.

For a smooth function $f:I\rightarrow {\mathbb R}$ defined on an open interval, we will also say that $f$ is
{\it strictly convex} if the graph of $f$ has  positive curvature $\kappa$ with respect to the upward unit normal $N$. This condition is equivalent to the positivity of $f''(x)$ on $I$.
\vskip 0.3cm

Suppose that $X$ is a strictly convex  curve  in the plane
 ${\mathbb R}^{2}$ with the unit normal $N$ pointing to the convex side.
 For a fixed point $P=A \in X$, and for a sufficiently small $h>0$, consider the  line $m$ passing through
 $P+hN(P)$ which is parallel to
 the tangent $\ell$ of $X$ at $P$.
 Let us denote by $A_1$ and $A_2$ the points where the line $m$ intersects the curve $X$.
We denote by $L_P(h)$ the length $|A_1A_2|$ of the chord $A_1A_2$.

Let us denote by  $\ell_1$, $\ell_2$ the tangent
lines passing through the points $A_1,A_2$ and by $B,B_1,B_2$
the intersection points $\ell_1\cap \ell_2$, $\ell\cap \ell_1$,
$\ell\cap \ell_2$, respectively. We denote by $T_P(h), U_P(h)$ the area  $|\bigtriangleup AA_1A_2|$ ,
$|\bigtriangleup BB_1B_2|$, of triangles, respectively.
 Then, obviously  we have $T_P(h)=\frac{h}{2}L_P(h)$.
\vskip 0.3cm

In this paper, first of all, in Section 2 we prove the following:
 \vskip 0.3cm

 \noindent {\bf Theorem 1.}
 Let $X$ denote a strictly convex $C^{(3)}$ curve  in the plane
 ${\mathbb R}^{2}$. Then we have

 \begin{equation}\tag{1.1}
   \begin{aligned}
 \lim_{h\rightarrow 0} \frac{T_P(h)}{h\sqrt{h}}=\frac{\sqrt{2}}{\sqrt{\kappa(P)}}
         \end{aligned}
   \end{equation}
 and
 \begin{equation}\tag{1.2}
   \begin{aligned}
 \lim_{h\rightarrow 0} \frac{U_P(h)}{h\sqrt{h}}=\frac{\sqrt{2}}{2\sqrt{\kappa(P)}}.
         \end{aligned}
   \end{equation}

 \vskip 0.3cm
Next in Section 3, using Theorem 1 we characterize parabolas as follows.
 \vskip 0.3cm
 \noindent {\bf Theorem 2.}
 Let $X=X(s)$ denote a strictly convex $C^{(3)}$ curve  in the plane
 ${\mathbb R}^{2}$. Suppose that for all $s$ and sufficiently small $h_i, i=1,2$,
 the curve $X$ satisfies
 \begin{equation}\tag{1.3}
   \begin{aligned}
 U(s,h_1,h_2)=\lambda(s)T(s,h_1,h_2).
         \end{aligned}
   \end{equation}
  Then,  we have $\lambda(s)=\frac{1}{2}$ and $X$ is an open part of a parabola.
\vskip 0.3cm
In [6], Krawczyk showed that for a strictly convex $C^{(4)}$ curve $X=X(s)$ in the plane
 ${\mathbb R}^{2}$, the following holds:
  \begin{equation}\tag{1.4}
   \begin{aligned}
 \lim_{h_1,h_2 \rightarrow 0} \frac{T(s,h_1,h_2)}{U(s,h_1,h_2)}=2.
         \end{aligned}
   \end{equation}
His application of (1.4) states  that if a strictly convex $C^{(4)}$ curve $X=X(s)$ in the plane
 ${\mathbb R}^{2}$ satisfies (1.3), then $\lambda(s)=\frac{1}{2}$ and $X$ is an open part of
 the graph of a quadratic polynomial.

But, for example, consider   a function $f(x)$ given  by
 \begin{equation}\tag{3.14}
  \begin{aligned}
y=\frac{2\sqrt{a}cx+1-\sqrt{4\sqrt{a}cx+1}}{2c^2},
   \end{aligned}
  \end{equation}
where $a, c>0$. Then, the function $f$ is defined on $I=(-\frac{1}{4\sqrt{a}c}, \infty)$.
Its graph $X$ is strictly convex and  satisfies (1.3) with $\lambda= \frac{1}{2}$.
Note that $X$ is not
the graph of a quadratic polynomial, but an open part of the parabola given in (3.16) in Section 3.
\vskip 0.30cm
In [5], the first
author and Y. H. Kim established five characterizations of parabolas, which are the converses of
well-known properties of parabolas  originally due to Archimedes
([8]).  In [3] and [4], they  also proved  the  higher dimensional analogues of some results in [5].
\vskip 0.30cm
Among  the  graphs of functions,  B.  Richmond and T. Richmond  established a dozen characterizations of parabolas
using elementary techniques ([7]). In their paper, parabola means the graph of a quadratic polynomial in one
variable.
 \vskip 0.30cm
Finally, we pose a question  as follows.
\vskip 0.30cm
\noindent {\bf Question 3.}
Let X be a strictly convex $C^{(3)}$ plane curve.
Suppose that for each $P\in X$ there exists a positive number $\epsilon=\epsilon(P)>0$
such that $U_P(h)=T_P(h)/2$ for all $h$ with $0<h<\epsilon(P)$.
Then, is it an open part of a parabola?

\vskip 0.3cm
  Throughout this article, all curves are of class $C^{(3)}$ and connected, unless otherwise mentioned.
  \vskip 0.50cm

 \section{Preliminaries and  Theorem 1}
 \vskip 0.50cm

In order to prove Theorem 1, we  need the following lemma ([5]).
\vskip 0.3cm

 \noindent {\bf Lemma 4.} Suppose that   $X$  is a strictly convex  $C^{(3)}$ curve  in the plane
 ${\mathbb R}^{2}$ with  the unit normal $N$ pointing to the convex side.  Then
 we have
  \begin{equation}\tag{2.1}
   \begin{aligned}
   \lim_{h\rightarrow 0} \frac{1}{\sqrt{h}}L_P(h)= \frac{2\sqrt{2}}{\sqrt{\kappa(P)}},
    \end{aligned}
   \end{equation}
 where $\kappa(P)$ is the curvature of $X$ at $P$ with respect to  the unit normal $N$.
 \vskip 0.3cm

First of all, we give a proof of (1.1) in  Theorem 1. Since $T_P(h)=\frac{h}{2}L_P(h)$, it
follows from Lemma 4 that the following holds:
 \begin{equation}\tag{1.1}
   \begin{aligned}
   \lim_{h\rightarrow 0} \frac{1}{h\sqrt{h}}T_P(h)= \frac{\sqrt{2}}{\sqrt{\kappa(P)}}.
    \end{aligned}
   \end{equation}
\vskip 0.3cm
In order to prove (1.2) in Theorem 1, we  fix an arbitrary  point $P$ on $X$.
Then, we may take a coordinate system $(x,y)$
 of  ${\mathbb R}^{2}$: $P$ is taken to be the origin $(0,0)$ and $x$-axis is the tangent line $\ell$ of $X$ at $P$.
 Furthermore, we may regard $X$ to be locally  the graph of a non-negative strictly convex  function $f: {\mathbb R}\rightarrow {\mathbb R}$
 with $f(0)=f'(0)=0$. Then $N$ is the upward unit normal.

 Since the curve $X$ is of class $C^{(3)}$, the Taylor's formula of $f(x)$ is given by
 \begin{equation}\tag{2.2}
 f(x)= ax^2 + g(x),
  \end{equation}
where  $2a=f''(0)$ and $g(x)$ is an $O(|x|^3)$  function.
From  $\kappa(P)=f''(0)>0$, we see that $a$ is positive.

For a sufficiently small $h>0$, we denote by $A_1(s,f(s))$ and $A_2(t,f(t))$
 the points where the line $m:y=h$ meets the curve $X$ with $s<0<t$.
 Then $f(s)=f(t)=h$ and we get $B_1(s-h/f'(s),0)$, $B_2(t-h/f'(t),0)$ and $B(x_0,y_0)$ with
  \begin{equation}\tag{2.3}
x_0=\frac{tf'(t)-sf'(s)}{ f'(t)-f'(s)}
  \end{equation}
and
  \begin{equation}\tag{2.4}
y_0=\frac{(t-s)f'(t)f'(s)+h(f'(t)-f'(s))}{ f'(t)-f'(s)}.
  \end{equation}
Noting that $L_P(h)=t-s$, one obtains
 \begin{equation}\tag{2.5}
  \begin{aligned}
2U_P(h)&=\{t-s -\frac{h}{f'(t)}+\frac{h}{f'(s)}\}(-y_0)\\
&=h^2\frac{(f'(t)-f'(s))}{-f'(s)f'(t)}-2hL_P(h)+\frac{-f'(s)f'(t)}{ f'(t)-f'(s)}L_P(h)^2.
   \end{aligned}
  \end{equation}

  It follows from (2.5) that
\begin{equation}\tag{2.6}
  \begin{aligned}
2\frac{U_P(h)}{h\sqrt{h}}&=\alpha_P(h)-2\frac{L_P(h)}{\sqrt{h}}+ \frac{1}{\alpha_P(h)}(\frac{L_P(h)}{\sqrt{h}})^2,
   \end{aligned}
  \end{equation}
  where we denote
 \begin{equation}\tag{2.7}
  \begin{aligned}
\alpha_P(h)=\sqrt{h}\frac{(f'(t)-f'(s))}{-f'(s)f'(t)}.
   \end{aligned}
  \end{equation}

Finally, we prove a lemma, which together with (2.6) and Lemma 4, completes the proof of (2) in Theorem 1.
\vskip 0.3cm
   \noindent {\bf Lemma 5.} We have the following.
 \begin{equation}\tag{2.8}
  \begin{aligned}
\lim_{h\rightarrow0}\alpha_P(h)=\frac{\sqrt{2}}{\sqrt{\kappa(P)}}.
   \end{aligned}
  \end{equation}

  \noindent {\bf Proof.}
  Note that
  \begin{equation}\tag{2.9}
  \begin{aligned}
\alpha_P(h)=\frac{\beta_P(h)}{\gamma_P(h)},
   \end{aligned}
  \end{equation}
where we denote
    \begin{equation}\tag{2.10}
  \begin{aligned}
\beta_P(h)=\frac{f'(t)-f'(s)}{t-s}
   \end{aligned}
  \end{equation}
and
  \begin{equation}\tag{2.11}
  \begin{aligned}
\gamma_P(h)=\frac{-f'(s)f'(t)}{\sqrt{h}(t-s)}.
   \end{aligned}
  \end{equation}

 Applying mean value theorem to the derivative $f'(x)$ of $f(x)$ shows 
 that as $h$ tends to $0$, $\beta_P(h)$ goes to $f''(0)=\kappa(P)$. 
 To get the limit of $\gamma_P(h)$,
  we put
   \begin{equation}\tag{2.12}
  \begin{aligned}
\delta_P(h)=\frac{f'(s)f'(t)}{st}
   \end{aligned}
  \end{equation}
  and
   \begin{equation}\tag{2.13}
  \begin{aligned}
\eta_P(h)=\frac{-st}{\sqrt{h}(t-s)}.
   \end{aligned}
  \end{equation}
  Then, we have
     \begin{equation}\tag{2.14}
  \begin{aligned}
\gamma_P(h)=\delta_P(h)\eta_P(h).
   \end{aligned}
  \end{equation}

 Note that
    \begin{equation}\tag{2.15}
  \begin{aligned}
\lim_{h\rightarrow 0}\delta_P(h)=\kappa(P)^2.
   \end{aligned}
  \end{equation}
 If we use $L_P(h)=t-s$, $\eta_P(h)$ can be written as
    \begin{equation}\tag{2.16}
  \begin{aligned}
\eta_P(h)=(-\frac{st}{h})/(\frac{L_P(h)}{\sqrt{h}})
   \end{aligned}
  \end{equation}
and the numerator of (2.16) can be decomposed as
   \begin{equation}\tag{2.17}
  \begin{aligned}
-\frac{st}{h}=(\frac{L_P(h)}{\sqrt{h}}-\frac{t}{\sqrt{h}})\frac{t}{\sqrt{h}}.
   \end{aligned}
  \end{equation}

  Now, to obtain the limit of $\frac{t}{\sqrt{h}}$, we use (2.2).
   Recalling that $\kappa(P)=f''(0)=2a$, we have
   \begin{equation}\tag{2.18}
  \begin{aligned}
\frac{t}{\sqrt{h}}=\frac{t}{\sqrt{at^2+g(t)}}.
   \end{aligned}
  \end{equation}
 Since $g(x)$ is an $O(|x|^3)$  function,  (2.18) implies that
  $ \lim_{h\rightarrow 0}t/\sqrt{h}=1/\sqrt{a}$.
 Hence, together with (2.17), Lemma 4 shows that $ \lim_{h\rightarrow 0}(-st)/h=1/a$,
  and hence from (2.16) we get
  \begin{equation}\tag{2.19}
  \begin{aligned}
  \lim_{h\rightarrow 0}\eta_P(h)=\frac{1}{2\sqrt{a}}.
   \end{aligned}
  \end{equation}
Thus, it follows from (2.14) and (2.15) that
  \begin{equation}\tag{2.20}
  \begin{aligned}
  \lim_{h\rightarrow 0}\gamma_P(h)=2a\sqrt{a}.
   \end{aligned}
  \end{equation}
  
  Using $\kappa(P)=2a$, together with (2.9) and (2.10), (2.20) completes the proof of Lemma 5. $\quad \square$
 \vskip 0.50cm

\section{Proof of Theorem 2}
 \vskip 0.3cm
In this section, we prove Theorem 2.

 Suppose that $X=X(s)$ denote a strictly convex $C^{(3)}$ curve  in the plane
 ${\mathbb R}^{2}$ which  satisfies for all $s$ and sufficiently small $h_i, i=1,2$,
  \begin{equation}\tag{1.3}
   \begin{aligned}
 U(s,h_1,h_2)=\lambda(s)T(s,h_1,h_2).
         \end{aligned}
   \end{equation}
  Then, in particular, for all $P=X(s)$ and sufficiently small $h>0$ the curve $X$ satisfies
   \begin{equation}\tag{3.1}
   \begin{aligned}
 U_P(h)=\lambda(P)T_P(h).
         \end{aligned}
   \end{equation}
  Hence, Theorem 1 implies that $\lambda(P)=\frac{1}{2}$.
 \vskip 0.3cm
In order to prove the remaining part of Theorem 2,
first, we  fix an arbitrary  point $A$ on $X$.
As in Section 1,  we take a coordinate system $(x,y)$
 of  ${\mathbb R}^{2}$: $A$ is taken to be the origin $(0,0)$
 and $x$-axis is the tangent line $\ell$ of $X$ at $A$.
 Furthermore, we may regard $X$ to be locally  the graph of a non-negative strictly convex
  function $f: {\mathbb R}\rightarrow {\mathbb R}$
 with $f(0)=f'(0)=0$ and $2a=f''(0)>0$.

For sufficiently small $|s|$ and $|t|$ with $0<s<t$ or $t<s<0$,
we let $A_1=(s,f(s)),A_2=(t,f(t))$ be two  neighboring points of $A$ on $X$.
Then, the area $T(s,t)$ of the triangle $\bigtriangleup AA_1A_2$ is given by
 \begin{equation}\tag{3.2}
  \begin{aligned}
2\epsilon T(s,t)&=(sf(t)-tf(s)),
   \end{aligned}
  \end{equation}
where $\epsilon=1$ if $0<s<t$ and $\epsilon=-1$ if $t<s<0$.

Denote by $\ell$, $\ell_1$, $\ell_2$ the tangent
lines passing through the points $A,A_1,A_2$ and by $B,B_1,B_2$
the intersection points $\ell_1\cap \ell_2$, $\ell\cap \ell_1$,
$\ell\cap \ell_2$, respectively.
Then we have
$B_1(s-f(s)/f'(s),0)$, $B_2(t-f(t)/f'(t),0)$ and $B(x_0,y_0)$ with
   \begin{equation}\tag{3.3}
y_0=\frac{(t-s)f'(t)f'(s)+f(s)f'(t)-f'(s)f(t)}{ f'(t)-f'(s)}.
  \end{equation}
Hence the area $U(s,t)$ of the triangle $\bigtriangleup BB_1B_2$ is given by
 \begin{equation}\tag{3.4}
  \begin{aligned}
2\epsilon U(s,t)&=\{t-s -\frac{f(t)}{f'(t)}+\frac{f(s)}{f'(s)}\}(-y_0)\\
&=\frac{\{(t-s)f'(t)f'(s)+f(s)f'(t)-f'(s)f(t)\}^2}{ f'(s)f'(t)(f'(t)-f'(s))}.
   \end{aligned}
  \end{equation}

Second, we prove
\vskip 0.3cm
\noindent {\bf Lemma 6.} The function $f$ satisfies the following:
\begin{equation}\tag{3.5}
  \begin{aligned}
f(t)f'(t)^2=4a(tf'(t)-f(t))^2,
   \end{aligned}
  \end{equation}
where  $a$ is given by $f''(0)=2a$.

\vskip 0.3cm
\noindent {\bf Proof.} Since the curve $X$ satisfies (1.3) with $\lambda=1/2$, we get
$2U(s,t)=T(s,t)$. By letting $s\rightarrow 0$, from (3.2) we get
 \begin{equation}\tag{3.6}
  \begin{aligned}
\epsilon \lim_{s\rightarrow 0}\frac{T(s,t)}{s}=\frac{f(t)}{2},
   \end{aligned}
  \end{equation}
where we use $f'(0)=0$. From (3.4) we also get
 \begin{equation}\tag{3.7}
  \begin{aligned}
2\epsilon \lim_{s\rightarrow 0}\frac{U(s,t)}{s}&=\frac{\{f''(0)tf'(t)-f''(0)f(t)\}^2}{ f''(0)f'(t)^2}\\
&=2a\frac{\{tf'(t)-f(t)\}^2}{ f'(t)^2},
   \end{aligned}
  \end{equation}
where we use $f'(0)=0$ and $f''(0)=2a>0$. Together with (3.6), (3.7) completes the proof. $\quad \square$
\vskip 0.3cm

Third, we prove
\vskip 0.3cm
\noindent {\bf Lemma 7.} The function $f$ satisfies the following:
\begin{equation}\tag{3.8}
  \begin{aligned}
2f(t)^2f''(t)=f'(t)^2\{tf'(t)-f(t)\}.
   \end{aligned}
  \end{equation}

\vskip 0.3cm
\noindent {\bf Proof.} By letting $s\rightarrow t$,  we get from (3.2)
 \begin{equation}\tag{3.9}
  \begin{aligned}
\epsilon \lim_{s\rightarrow t}\frac{T(s,t)}{s-t}&=\frac{1}{2}\lim_{s\rightarrow t}\frac{sf(t)-tf(s)}{s-t}\\
&=\frac{1}{2}(f(t)-tf'(t)).
   \end{aligned}
  \end{equation}
On the other hand, from (3.4) we get
 \begin{equation}\tag{3.10}
  \begin{aligned}
2\epsilon \lim_{s\rightarrow t}\frac{U(s,t)}{s-t}
&=\lim_{s\rightarrow t}\frac{\{(t-s)f'(t)f'(s)+f(s)f'(t)-f'(s)f(t)\}^2}{ (s-t)f'(s)f'(t)(f'(t)-f'(s))}\\
&=-\frac{f(t)^2f''(t)}{f'(t)^2}.
   \end{aligned}
  \end{equation}
Since $T=2U$, together with (3.9), (3.10) completes the proof.  $\quad \square$
\vskip 0.3cm

By eliminating $tf'(t)-f(t)$ from (3.5) and (3.8), we get
\begin{equation}\tag{3.11}
  \begin{aligned}
f''(t)=\frac{1}{4\sqrt{a}}\frac{f'(t)^3}{f(t)^{3/2}}.
   \end{aligned}
  \end{equation}
Letting $y=f(t)$,  a standard method of ordinary differential equations 
using the substitution $w=dy/dt$ and $y''(t)=w(dw/dy)$ leads to
\begin{equation}\tag{3.12}
  \begin{aligned}
dt=(\frac{1}{2\sqrt{ay}}+c)dy,
   \end{aligned}
  \end{equation}
where $c$ is a constant. Since $f(0)=0$, we obtain from (3.12)
\begin{equation}\tag{3.13}
  \begin{aligned}
t=\frac{1}{\sqrt{a}}(\sqrt{y}+cy).
   \end{aligned}
  \end{equation}

After replacing $t$ by $x$, we have for $y=f(x)$
\begin{equation}\tag{3.14}
 \begin{aligned}
 y= \begin{cases}
  \frac{2\sqrt{a}cx+1-\sqrt{4\sqrt{a}cx+1}}{2c^2}, & \text{if $c\ne0,$} \\
 ax^2, & \text{if $c= 0.$}
  \end{cases}
  \end{aligned}
  \end{equation}
  Note that
  \begin{equation}\tag{3.15}
  \begin{aligned}
f(0)=f'(0)=0 , f''(0)=2a \quad \text{and} \quad  f'''(0)=-12\sqrt{a}ac \quad \text{or} \quad 0.
   \end{aligned}
  \end{equation}
  It follows from (3.14) that  the curve $X$ around an arbitrary point $A$
  is an open part of  the parabola defined by
\begin{equation}\tag{3.16}
  \begin{aligned}
ax^2-2\sqrt{a}cxy+c^2y^2-y=0.
   \end{aligned}
  \end{equation}

 \vskip 0.3cm

 Finally using (3.15), in the same manner as in [5], we can show that the curve $X$ is globally an open
 part of a  parabola. This completes the proof of Theorem 2.

 \vskip 0.50cm

\section{Corollaries and examples}
 \vskip 0.5cm
In this section, we give some corollaries and examples.

Suppose that $X=X(s)$ is a strictly convex $C^{(3)}$ curve  in the plane
 ${\mathbb R}^{2}$ which  satisfies for all $s$ and sufficiently small $h_i, i=1,2$,
  \begin{equation}\tag{4.1}
   \begin{aligned}
 U(s,h_1,h_2)=\lambda(s)T(s,h_1,h_2)^{\mu(s)},
         \end{aligned}
   \end{equation}
   where $\lambda(s)$ and $\mu(s)$ are some functions.
     Then, in particular, for all $P=X(s)$ and sufficiently small $h>0$ the curve $X$ satisfies
   \begin{equation}\tag{4.2}
   \begin{aligned}
 U_P(h)=\lambda(P)T_P(h)^{\mu(P)}.
         \end{aligned}
   \end{equation}
 Using Theorem 1, by letting $h\rightarrow 0$ we see that  $\mu(P)=1$.
  Hence, Theorem 1 again implies that $\lambda(P)=\frac{1}{2}$.

  Thus, from Theorem 2 we get
 \vskip 0.3cm

\noindent {\bf Corollary 8.} Let $X$ denote   a strictly convex curve  in the plane ${\mathbb R}^{2}$.
Then, the following are equivalent.
\vskip 0.3cm
\noindent 1) $X$ satisfies (4.1) for some functions $\lambda(s)$ and $\mu(s)$.

\noindent 2) $X$ satisfies (4.1) with $\lambda=\frac{1}{2}$ and $\mu=1$ .

\noindent 3) $X$ is an open part of a parabola.

\vskip 0.3cm

 Finally, we give an example of a convex curve which satisfies
  \begin{equation}\tag{1.5}
   \begin{aligned}
U_P(h)=\frac{1}{2} T_P(h).
         \end{aligned}
   \end{equation}
  for  sufficiently small $h>0$ at every point $P\in X$,  but it is not a  parabola.
Note that the example is not of class $C^{(2)}$, and hence it is not strictly convex either.
\vskip 0.3cm

 \noindent {\bf Example 9.}
Consider the graph $X$ of a function $f:{\mathbb R}\rightarrow {\mathbb R}$ which is given by
for some positive distinct constants $a$ and $b$
\begin{equation}\tag{4.3}
 \begin{aligned}
 f(x)= \begin{cases}
  ax^2, & \text{if $x<0,$} \\
 bx^2, & \text{if $x\ge 0.$}
  \end{cases}
  \end{aligned}
  \end{equation}
 It is straightforward to show that if $P$ is the origin, then for all $h$ we have
\begin{equation}\tag{4.4}
   \begin{aligned}
 U_P(h)=\frac{1}{2}T_P(h).
         \end{aligned}
   \end{equation}
Hence $X$ satisfies $U_P(h)=T_P(h)/2$ at the origin  for all $h>0$. If $P\in X$ is not the origin, then
there  exists a positive number $\varepsilon(P)$ such that
 for every positive number $h$ with $h<\varepsilon(P)$,
$X$ satisfies $U_P(h)=T_P(h)/2$.

Thus, $X$ satisfies $U_P(h)=T_P(h)/2$ for sufficiently small $h>0$ at every point $P\in X$. But it is not a parabola.
\vskip 0.3cm

 \vskip 0.50cm


\vskip 1.0 cm

Department of Mathematics, \par Chonnam National University,\par
Kwangju 500-757, Korea

{\tt E-mail:dosokim@chonnam.ac.kr and mathtsim@naver.com} \vskip 0.3 cm


\vskip 0.50cm

\begin{thebibliography}{5.4}
\bibitem {}
Day, W. A., {\it
Inequalities for areas associated with conics},
Amer. Math. Monthly 98 (1991), no. 1, 36-39.


\bibitem {}
do Carmo, M. P., {\it Differential Geometry of
Curves and Surfaces}, Prentice-Hall, Englewood Cliffs,
NJ, 1976.

\bibitem {}
Kim, D.-S.  and  Kim, Y. H., {\it Some characterizations of spheres and elliptic paraboloids}, Linear Algebra Appl.
437 (2012), no. 1, 113-120.


\bibitem {}
Kim, D.-S.  and  Kim, Y. H., {\it Some characterizations of spheres and elliptic paraboloids II},
 Linear Algebra Appl., 438 (2013), no. 3, 1356-1364.

 \bibitem {}
Kim, D.-S.  and  Kim, Y. H., {\it On the Archimedean characterization of parabolas},  Bull. Korean Math. Soc.,
To appear. arXiv:1305.3337.


\bibitem {}
Krawczyk, J., {\it On areas associated with a curve},
Zesz. Nauk. Uniw.Opol. Mat. 29 (1995), 97-101.

\bibitem {}
Richmond, B. and Richmond, T.,  {\it How to recognize a parabola}, Amer. Math. Monthly 116(2009), no.10, 910-922.

\bibitem {}
Stein, S., {\it Archimedes. What did he do besides cry Eureka?},
Mathematical Association of America, Washington, DC, 1999.

\end{thebibliography}
\end{document}